# A note on the values of an Euler sum at negative integers and relation to a convolution of Bernoulli numbers

Khristo N. Boyadzhiev,
Ohio Northern University
Department of Mathematics
Ada, Ohio, 45810, USA
E-mail: k-boyadzhiev@onu.edu

H. Gopalkrishna Gadiyar and R. Padma
AU-KBC Research Centre
M. I. T. Campus of Anna University
Chromepet, Chennai 600 044, INDIA
E-mail: {gadiyar, padma}@au-kbc.org

**1. Introduction**

We consider the function

$$h(s) = \sum_{n=1}^{\infty} H_n \frac{1}{n^s}, \qquad (1.1)$$

where

$$H_n = 1 + \frac{1}{2} + \frac{1}{3} + \ldots + \frac{1}{n}, \qquad (1.2)$$

are the Harmonic numbers and $Re\, s > 1$. This function was studied by Apostol-Vu [1] and Matsuoka [5] who provided an analytic extension for all complex values of $s$ and discussed its values and poles at the negative integers. In this note we shall find a relation between the values $h(1-n)$ and the numbers $A_n$, $n = 1, 2, \ldots$ defined as the convolution.

$$A_n = \sum_{k+j=n} \frac{B(k)}{k!k} \frac{B(j)}{j!}, \quad k = 1, 2, \ldots; j = 0, 1, \ldots \qquad (1.3)$$

where $B(m) = B_m$ are the Bernoulli numbers for $m \neq 1$ and $B(1) = -B_1 = 1/2$. Thus

$$\frac{ze^z}{e^z-1} = \frac{-z}{e^{-z}-1} = \sum_{n=0}^{\infty} \frac{B(n)}{n!} z^n. \tag{1.4}$$

We have also

$$\log\left(\frac{e^z-1}{z}\right) = \sum_{n=1}^{\infty} \frac{B(n)}{n!\,n} z^n. \tag{1.5}$$

The product of the functions in (1.4) and (1.5) is the generating function of the numbers $A_n$,

$$\frac{ze^z}{e^z-1} \log\left(\frac{e^z-1}{z}\right) = \sum_{n=1}^{\infty} A_n z^n. \tag{1.6}$$

The relation between $A_n$ and the values of $h(s)$ is based on the evaluation of the following integral

$$F(s) = \frac{\Gamma(1-s)}{2\pi i} \int_L \frac{z^{s-1} e^z}{e^z-1} \log\left(\frac{e^z-1}{z}\right) dz, \tag{1.7}$$

where $L$ is the Hankel contour consisting of three parts: $L = L_- \cup L_+ \cup L_\epsilon$, with $L_-$ the "lower side" (i.e. $arg(z) = -\pi$) of the ray $(-\infty, -\epsilon)$, $\epsilon > 0$, traced left to right, and $L_+$ the "upper side" ($arg(z) = \pi$) of this ray traced right to left. Finally, $L_\epsilon = \{z = \epsilon e^{\theta i} : -\pi \leq \theta \leq \pi\}$ is a small circle traced counterclockwise and connecting the two sides of the ray. Such contours were used, for instance, in [2], [3].

We note that convolutions like (1.3) appear in the Matiyasevich version of Miki's identity - see [6].

## 2. Main results

The main results of this article are given in the folowingTheorem and the three corollaries

**Theorem.** For $Re\, s > 1$,

$$F(s) = h(s) - \zeta(s+1) + \psi(s)\zeta(s) + \zeta'(s), \tag{2.1}$$

where $\zeta(s)$ is the Riemann zeta function and $\psi(s) = \Gamma'(s)/\Gamma(s)$ is the Polygamma function.

As $F(s)/\Gamma(1-s)$ is an entire function (from (1.7), this provides an extension of the right hand side in (2.1) to all complex $s$.

(The proof of the theorem is given in Section 3.)

It is easy to see that when $s$ is a negative integer or zero, the integration in (1.7) can be reduced to $L_\epsilon$ only, as the integrals on $L_+$ and $L_-$ cancel each-other. This way for the coefficients $A_n$ of the Taylor series (4) we have for $n = 1, 2, \ldots$

$$(n-1)! A_n = F(1-n). \tag{2.2}$$

We shall evaluate the right hand side of (8) for $s = 1 - n$ by considering three cases, $n > 1$ odd, $n = 1$, and $n$ even. The results are organized in three corollaries. Before listing these corollaries, we recall two properties of the Riemann zeta function. For $m = 1, 2, \ldots,$

$$\zeta(-2m) = 0 \text{ and } \zeta(1-2m) = -\frac{B_{2m}}{2m}. \tag{2.3}$$

**Corollary 1**. Let first $n = 2m + 1$, $m > 0$. Then

$$(2m)! A_{2m+1} = h(-2m) - \zeta(1-2m) = \frac{1}{2} B_{2m}(1 + \frac{1}{2m}). \tag{2.4}$$

Proof. From (2.2) we have $(2m)! A_{2m+1} = F(-2m)$. In order to evaluate $F(-2m)$ we use the well-known property of the Polygamma function

$$\psi(s) = \psi(1-s) - \pi \cot \pi s \tag{2.5}$$

to write

$$\psi(s)\zeta(s) = \psi(1-s)\zeta(s) - \zeta(s) \pi \cot \pi s. \tag{2.6}$$

Now, for $s = -2m$ we have $\psi(1+2m)\zeta(-2m) = 0$ and

$$\zeta(s) \pi \cot \pi s \big|_{s=-2m} = \zeta'(-2m) \tag{2.7}$$

as follows from the Taylor expansion centered at $s = -2m$

$$\zeta(s) \pi \cot \pi s = \zeta'(-2m) + \frac{1}{2}\zeta''(-2m)(s+2m) + O((s+2m)^2). \tag{2.8}$$

Thus from (2.1) we find

$$F(-2m) = h(-2m) - \zeta(1-2m). \tag{2.9}$$

The values $h(-2m)$ were computed by Matsuoka [5] as

$$h(-2m) = -\frac{B_{2m}}{4m} + \frac{B_{2m}}{2}; \tag{2.10}$$

(note that Matsuoka worked with the function $f(s) = h(s) - \zeta(s+1)$). Therefore, equation (2.4) follows from (2.9) and (2.10).

We note that $h(-2m)$ was also evaluated in [1], but incompletely (missing the second term on the right hand side in (2.10).

Now consider the case $s = 0$ in (2.1), i.e. $n = 1$ in (2.2).

**Corollary 2**. In a neighborhood of zero,

$$h(s) = \frac{1}{2s} + \frac{1}{2}(1+\gamma) + O(s), \tag{2.11}$$

where $\gamma = -\psi(1)$ is the Euler constant.

*Proof.* As found in [1] and [5], the function $h(s)$ has a simple pole at $s = 0$ with residue ½. In order to establish (2.11) we need to evaluate $h(s) - 1/(2s)$ at zero. The functions $\zeta(s+1)$ and $\psi(s)\zeta(s)$ have at zero residues 1 and ½ correspondingly, and so the function

$$\zeta(s+1) - \psi(s)\zeta(s) - \frac{1}{2s} \tag{2.12}$$

does not have a pole at $s = 0$. Moreover, one easily finds that at zero

$$\zeta(s+1) - \psi(s)\zeta(s) - \frac{1}{2s} = \frac{\gamma}{2} + \zeta'(0) + O(s). \tag{2,13}$$

Next we rewrite (2.1) in the form

$$h(s) - \frac{1}{2s} = F(s) + (\zeta(s+1) - \psi(s)\zeta(s) - \frac{1}{2s}) - \zeta'(s) \tag{2.14}$$

and also compute the coefficient $A_1 = 1/2 = F(0)$ from (3). From (2.13) and (2,14) we find

$$(h(s) - \frac{1}{2s})|_{s=0} = \frac{1}{2}(1+\gamma), \tag{2.15}$$

which proves (2.11).

Finally, we compute $F(1-n)$ for even $n = 2m$.

**Corollary 3.** For $m = 2, 3, \ldots$, in a neighborhood of $s = 1 - 2m$ the function $h(s)$ is represented as

$$h(s) = \frac{\zeta(1-2m)}{s+2m-1} + (2m-1)!A_{2m} - \psi(2m)\zeta(1-2m) + O(s+2m-1), \qquad (2.16)$$

and in a neighborhood of $s = -1$,

$$h(s) = \frac{-1}{12(s+1)} - \frac{1}{8} + \frac{\gamma}{12} + O(s+1). \qquad (2.17)$$

*Proof.* Apostol-Vu [1] and Matsuoka [5], showed that the function $h(s)$ has simple poles at the negative odd integers $s = 1 - 2m$ with residues $\zeta(1-2m)$. The same is true for the function $\zeta(s)\pi\cot\pi s$, as follows from the Taylor expansion at $s = 1 - 2m$,

$$\zeta(s)\pi\cot\pi s = \zeta(1-2m)\frac{1}{s+2m-1} + \zeta'(1-2m) + O(s+2m-1). \qquad (2.18)$$

We substitute this in (2.1) written in the form

$$h(s) = \zeta(s)\pi\cot\pi s + F(s) + \zeta(s+1) - \psi(1-s)\zeta(s) - \zeta'(s), \qquad (2.19)$$

to get

$$h(s) - \frac{\zeta(1-2m)}{s+2m-1} = \qquad (2.20)$$

$$F(s) + \zeta(s+1) - \psi(1-s)\zeta(s) - \zeta'(s) + \zeta'(1-2m) + O(s+2m-1).$$

From here, evaluating both sides for $s = 1 - 2m$,

$$[h(s) - \frac{\zeta(1-2m)}{s+2m-1}]|_{s=1-2m} = \qquad (2.21)$$

$$F(1-2m) + \zeta(2-2m) - \psi(2m)\zeta(1-2m) + O(s+2m-1)$$

and as $F(1-2m) = (2m-1)!A_{2m}$ and $\zeta(2-2m) = 0$, we obtain (2.16).

When $m = 1$, we have $\zeta(2-2m) = \zeta(0) = -1/2$, $\zeta(-1) = -1/12$, $\psi(2) = 1 - \gamma$, and by

direct computation from (1.3), $A_2 = 7/24$. Thus (24) follows from (2.21).

### 3. Proof of the theorem

Here we evaluate the integral in (1.7)

$$I(s) = \frac{1}{2\pi i} \int_L \frac{z^{s-1} e^z}{e^z - 1} \log\left(\frac{e^z - 1}{z}\right) dz, \tag{3.1}$$

where the contour $L$ is described in the introduction. We choose $\operatorname{Re} s > 1$ and set $\epsilon \to 0$. The integral over $L_\epsilon$ becomes zero, as the function

$$\frac{z e^z}{e^z - 1} \log\left(\frac{e^z - 1}{z}\right) \tag{3.2}$$

is holomorphic in a neighborhood of zero. Noticing that $z = x e^{-\pi i}$ on $L_-$ and $z = x e^{\pi i}$ on $L_+$, we find that

$$-I(s) = \frac{e^{-\pi i}}{2\pi i} \int_\infty^0 \frac{x^{s-1} e^{-x}}{1 - e^{-x}} \log\left(\frac{1 - e^{-x}}{x}\right) dx + \frac{e^{\pi i}}{2\pi i} \int_0^\infty \frac{x^{s-1} e^{-x}}{1 - e^{-x}} \log\left(\frac{1 - e^{-x}}{x}\right) dx$$

$$= \frac{\sin \pi s}{\pi} \int_0^\infty \frac{x^{s-1} e^{-x}}{1 - e^{-x}} \log\left(\frac{1 - e^{-x}}{x}\right) dx. \tag{3.3}$$

Next,

$$\int_0^\infty \frac{x^{s-1} e^{-x}}{1 - e^{-x}} \log\left(\frac{1 - e^{-x}}{x}\right) dx = \int_0^\infty \frac{x^{s-1} e^{-x}}{1 - e^{-x}} \log(1 - e^{-x}) dx - \int_0^\infty \frac{x^{s-1}}{e^x - 1} \log x \, dx. \tag{3.4}$$

We shall evaluate the two integrals on the right hand side in (3.4) one by one. First we use the expansion

$$\frac{\log(1 - e^{-x})}{1 - e^{-x}} = -\sum_{n=1}^\infty H_n e^{-nx} \tag{3.5}$$

(see [4, (7.57), p. 352]). Multiplying this by $x^{s-1}e^{-x}$ and integrating from zero to infinity, we find that

$$\int_0^\infty \frac{x^{s-1}e^{-x}}{1-e^{-x}} \log(1-e^{-x}) dx = -\sum_{n=1}^\infty H_n \int_0^\infty x^{s-1} e^{-(n+1)x} dx$$

$$= -\Gamma(s) \sum_{n=1}^\infty \frac{H_n}{(n+1)^s} = -\Gamma(s)(h(s) - \zeta(s+1)). \tag{3.6}$$

Next, differentiating for $s$ the representation

$$\Gamma(s)\zeta(s) = \int_0^\infty \frac{x^{s-1}}{e^x - 1} dx, \tag{3.7}$$

we obtain

$$\int_0^\infty \frac{x^{s-1}}{e^x - 1} \log x \, dx = \Gamma'(s)\zeta(s) + \Gamma(s)\zeta'(s) = \Gamma(s)(\psi(s)\zeta(s) + \zeta'(s)). \tag{3.8}$$

From (3.4), (3.6) and (3.8)

$$\int_0^\infty \frac{x^{s-1}e^{-x}}{1-e^{-x}} \log\left(\frac{1-e^{-x}}{x}\right) dx = -\Gamma(s)(h(s) - \zeta(s+1) + \psi(s)\zeta(s) + \zeta'(s)), \tag{3.9}$$

and therefore,

$$I(s) = \frac{1}{\pi} \Gamma(s) \sin(\pi s) (h(s) - \zeta(s+1) + \psi(s)\zeta(s) + \zeta'(s)). \tag{3.10}$$

Finally, (2.1) follows from here in view of the identity

$$\Gamma(s)\Gamma(1-s) = \frac{\pi}{\sin \pi s}. \tag{3.11}$$

The proof was done for $\operatorname{Re} s > 1$, but since the integral in (2.1) exists for all complex $s$, the representation hold for all $s$.